\documentclass[12pt,reqno]{article}

\usepackage{amsmath,amsthm,amsfonts,amssymb,latexsym,mathrsfs,color}

\newtheorem{thm}{Theorem}[section]
\newtheorem{lem}[thm]{Lemma}

\newtheorem{prop}[thm]{Proposition}
\newtheorem{conj}[thm]{Conjecture}

\newtheorem{pro}[thm]{Problem}

\theoremstyle{definition}

\newtheorem{ex}[thm]{Example}

\newtheorem{rem}[thm]{Remark}

\numberwithin{equation}{section}

\def\JS{\text{JS}}
\def\Jc{\text{Jc}}

\newcommand{\msn}{{\mathcal S}_n}
\newcommand{\exc}{{\rm exc\,}}

\newcommand{\rz}{{\rm RZ}}

\newcommand{\sep}{\preceq}

\newcommand{\as}{{\rm as\,}}

\title{$q$-log-convexity from linear transformations and polynomials with
only real zeros
\thanks{Supported partially by the National Natural Science Foundation of China (No.
11571150).
\newline\hspace*{5mm}
   {\it Email address:} bxzhu@jsnu.edu.cn (B.-X. Zhu)}}
\author{Bao-Xuan Zhu}
\date{\footnotesize School of Mathematics and Statistics,
         Jiangsu Normal University,
         Xuzhou 221116, PR China}

\begin{document}

\maketitle

\begin{abstract}
In this paper, we mainly study the stability of iterated polynomials
and linear transformations preserving the strong $q$-log-convexity
of polynomials

Let $[T_{n,k}]_{n,k\geq0}$ be an array of nonnegative numbers. We
give some criteria for the linear transformation
$$y_n(q)=\sum_{k=0}^nT_{n,k}x_k(q)$$
preserving the strong $q$-log-convexity (resp. log-convexity). As
applications, we derive that some linear transformations (for
instance, the Stirling transformations of two kinds, the
Jacobi-Stirling transformations of two kinds, the Legendre-Stirling
transformations of two kinds, the central factorial transformations,
and so on) preserve the strong $q$-log-convexity (resp.
log-convexity) in a unified manner. In particular, we confirm a
conjecture of Lin and Zeng, and extend some results of Chen {\it et
al.}, and Zhu for strong $q$-log-convexity of polynomials, and some
results of Liu and Wang for transformations preserving the
log-convexity.

The stability property of iterated polynomials implies the
$q$-log-convexity. By applying the method of interlacing of zeros,
we also present two criteria for the stability of the iterated Sturm
sequences and $q$-log-convexity of polynomials. As consequences, we
get the stabilities of iterated Eulerian polynomials of type $A$ and
$B$, and their $q$-analogs. In addition, we also prove that the
generating functions of alternating runs of type $A$ and $B$, the
longest alternating subsequence and up-down runs of permutations
form a $q$-log-convex sequence, respectively.
\end{abstract}

\section{Introduction}
The main objective of this paper is twofold: one is to study the
linear transformations preserving the strong $q$-log-convexity of
the sequences of polynomials and the other is to study one strong
property of $q$-log-convexity called the stability of polynomials.

 The
Jacobi-Stirling numbers $\JS_n^k( z)$ of the second kind, which were
introduced in \cite{EKLWY07}, are the coefficients of the integral
composite powers of the Jacobi differential operator
$$\ell_{\alpha,\beta}[y](t)=\frac{1}{(1-t)^\alpha(1+t)^\beta}\left( -(1-t)^{\alpha+1}(1+t)^{\beta+1} y'(t)\right)' ,$$
with fixed real parameters $\alpha,\beta\geq-1$. They also satisfy
the following recurrence relation:
\begin{equation*}
\left\{ \begin{array}{l} \JS_0^0(z)=1, \qquad \JS_n^k(z)=0, \quad  \text{if} \ k \not\in\{1,\ldots,n\}, \\
\JS_n^k(z)= \JS_{n-1}^{k-1}(z)+k(k+z)\,\JS_{n-1}^{k}(z),  \quad n,k
\geq 1,\end{array} \right.
 \end{equation*}
 where $z=\alpha+\beta+1.$
Similarly, the Jacobi-Stirling numbers $\Jc_n^k( z)$ of the first
kind are defined by
\begin{equation*}
\left\{ \begin{array}{l} \Jc_0^0(z)=1, \qquad \Jc_n^k(z)=0, \quad  \text{if} \ k \not\in\{1,\ldots,n\}, \\
\Jc_n^k(z)= \Jc_{n-1}^{k-1}(z)+(n-1)(n-1+z)\,\Jc_{n-1}^{k}(z), \quad
n,k \geq 1,\end{array} \right.
 \end{equation*}
 where $z=\alpha+\beta+1.$
 Actually, these numbers are a generalization of the Legendre-Stirling numbers of two kinds: it suffices to choose $\alpha
 =\beta=0$. Recently, the
Jacobi-Stirling numbers and Legendre-Stirling numbers have generated
a significant amount of interest from some researchers in
combinatorics, see~Andrews {\it et al.}~\cite{AEGL11,AL09}, Everitt
{\it et al.}~\cite{EKLWY07}, Mongelli~\cite{Monge12}, Lin and Zeng
\cite{LZ14} and Zhu \cite{Zhu14} for instance. In \cite{LZ14}, Lin
and Zeng proposed the next conjecture.
\begin{conj}\cite{LZ14}\label{Jacobi}
The Jacobi-Stirling transformations of two kinds
\begin{equation*}
y_n=\sum_{k=0}^{n}\JS_n^k(z)x_k\quad and\quad
w_n=\sum_{k=0}^{n}\Jc_n^k(z)x_k
 \end{equation*} preserve the log-convexity for $z=0,1$.
\end{conj}

Recall some notation and definitions. Let $\{a_n\}_{n\geq0}$ be a
sequence of nonnegative real numbers. It is called {\it log-convex}
(resp. {\it log-concave}) if for all $k\ge 1$, $a_{k-1}a_{k+1}\ge
a_k^2$ (resp. $a_{k-1}a_{k+1}\le a_k^2$), which is equivalent to
that $a_{n-1}a_{m+1}\ge a_na_m$ (resp. $a_{n-1}a_{m+1}\le a_na_m$)
for all $1\le n\le m$. The log-concave and log-convex sequences
arise often in combinatorics, algebra, geometry, analysis,
probability and statistics and have been extensively investigated,
see Stanley~\cite{Sta89}, Brenti~\cite{Bre94}, Liu and
Wang~\cite{LW07} and Zhu \cite{Zhu13} for details.

For two polynomials with real coefficients $f(q)$ and $g(q)$, denote
$f(q)\geq_q g(q)$ if the difference
$f\left(q\right)-g\left(q\right)$ has only nonnegative coefficients.
For a polynomial sequence $\{f_n(q)\}_{n\geq 0}$, it is called {\it
$q$-log-concave} first suggested by Stanley, if
$$f_n(q)^2- f_{n+1}(q)f_{n-1}(q)\geq_q0$$ for $n\geq 1$ and is called {\it strongly
$q$-log-concave} introduced by Sagan, if
$$ f_n(q)f_m(q)-f_{n+1}(q)f_{m-1}(q)\geq_q0$$
for any $n\geq m\geq1$. Obviously, the strong $q$-log-concavity
 of polynomials implies the
$q$-log-concavity. However, the converse does not hold.  The
$q$-log-concavity of polynomials have been extensively studied, see
Butler \cite{But90}, Leroux \cite{Le90}, and Sagan \cite{Sag921} for
instance.

For the polynomial sequence $\{f_n(q)\}_{n\geq 0}$, it is called
{\it $q$-log-convex} introduced by Liu and Wang, if
$$ f_{n+1}(q)f_{n-1}(q)-f_n(q)^2\geq_q 0$$ for $n\geq 1$ and is called {\it strongly $q$-log-convex}
defined by Chen {\it et al.}, if
$$f_{n+1}(q)f_{m-1}(q)-  f_n(q)f_m(q)\geq_q0$$
for any $n\geq m\geq1$. Clearly, strong $q$-log-convexity of
polynomials implies the $q$-log-convexity. However, the converse
does not hold.

The operator theory often is used to study the log-concavity or
log-convexity. For example, the log-convexity and log-concavity are
preserved under the binomial convolution respectively, see Davenport
and P\'olya~\cite{DP49} and Wang and Yeh~\cite{WY07}. Br\"and\'en
\cite{Bra06} studied some linear transformations preserving the
P\'olya frequency property of sequences. Brenti \cite{Bre89}
obtained some transformations preserving the log-concavity. Liu and
Wang \cite{LW07} also studied linear transformations preserving the
log-convexity. In \cite{ZS15,Zhu17}, we strengthened partial results
for the linear transformations preserving the log-convexity to the
strong $q$-log-convexity. However, there are fewer results about the
linear transformations preserving the strong $q$-log-convexity. One
of the aims of this paper is to continue studying linear
transformations preserving the strong $q$-log-convexity

 Given an array
$[T_{n,k}]_{n,k\geq0}$ of nonnegative real numbers and a sequence of
polynomials $\{x_n(q)\}_{n\geq0}$, define the polynomials
\begin{equation*}\label{q-p}
y_n(q)=\sum_{k\geq0}T_{n,k}x_k(q)
\end{equation*}
for $n\geq0$. If we take $x_k(q)=q^k$, then it was demonstrated that
the corresponding sequence $\{y_n(q)\}_{n\geq0}$ has the
$q$-log-convexity or strong $q$-log-convexity for many famous
triangles $[T_{n,k}]_{n,k\geq0}$, including the Stirling triangle of
the second kind, the Jacobi-Stirling triangle of the second kind,
the Legendre-Stirling triangle of the second kind, the Eulerian
triangles of type $A$ and $B$, the Narayana triangles of type $A$
and $B$, and so on, see
\cite{CTWY10,CWY10,CWY11,LW07,LZ15,Zhu13,Zhu14} for instance. Thus
it is natural to consider the strong $q$-log-convexity of the linear
transformation $y_n(q)$ by that of $x_n(q)$. On the other hand, note
that a log-convex sequence is one special case of the strongly
$q$-log-convex sequence since the real number sequence
$\{a_n\}_{n\geq0}$ is log-convex if and only if $a_{n-1}a_{m+1}\ge
a_na_m$ for all $1\le n\le m$. So it is easy to see that the linear
transformation preserving the strong $q$-log-convexity also
preserves the log-convexity.

Let $[T_{n,k}]_{n,k\geq0}$ be an array of nonnegative numbers
satisfying the recurrence relation:
\begin{equation}\label{rr+Chen}
T_{n,k}=(a_0n+a_2k+a_3)T_{n-1,k}+(b_0n+b_2k+b_3)T_{n-1,k-1}
\end{equation}
with $T_{n,k}=0$ unless $0\le k\le n$ and $T_{0,0}=1$. Chen {\it et
al.} \cite{CWY11} proved the strong $q$-log-convexity of the row
generating functions $T_n(q)=\sum_{k=0}^n T_{n,k}q^k$ when all
$a_0,a_2,b_0$ and $b_2$ are nonnegative real numbers.

In \cite{Zhu14}, we considered another array of nonnegative numbers
$[T_{n,k}]_{n,k\geq0}$ satisfying the recurrence relation:
\begin{equation}\label{rr+Zhu}
T_{n,k}=(a_1k^2+a_2k+a_3)T_{n-1,k}+(b_1k^2+b_2k+b_3)T_{n-1,k-1},
\end{equation}
where $T_{0,0}=1$ and $T_{0,k}=0$. We also proved the strong
$q$-log-convexity of its row generating functions when all
$a_0,a_1,b_0$ and $b_1$ are nonnegative real numbers.

In this paper, we consider a more generalized array of nonnegative
numbers $[T_{n,k}]_{n,k\geq0}$ satisfying the recurrence relation:
\begin{equation}\label{rr+Gen}
T_{n,k}=[r(n)+f(k)]\,T_{n-1,k}+[s(n)+g(k)]\,T_{n-1,k-1}+[t(n)+h(k)]\,T_{n-1,k-2}
\end{equation}
with $T_{n,k}=0$ unless $0\le k\le n$ and $T_{0,0}=1$. For $n\geq0$,
let $T_n(q)=\sum_{k=0}^n T_{n,k}q^k$ be the row generating
functions.

Recall that a matrix $M=(m_{ij})_{i,j\geq 0}$ of nonnegative numbers
is said to be {\it $r$-order totally positive} (TP$_r$ for short) if
its all minors of order at most $r$ are nonnegative.  Total
positivity of matrices has been extensively studied and is very
useful, see Karlin~\cite{Kar68} for more details. By total
positivity of matrices, we have the next extensive result for linear
transformations preserving the strong $q$-log-convexity.
\begin{thm}\label{thm+transf+squar}
Let $[T_{n,k}]_{n,k\ge 0}$ be the nonnegative array satisfying the
recurrence (\ref{rr+Gen}). Assume that the matrix $[T_{n,k}]_{n,k\ge
0}$ is TP$_2$ and all $r(n),s(n),t(n),f(n),g(n)$ and $h(n)$ are
nonnegative and increasing in $n$ for $n\geq0$. If
$\{x_n(q)\}_{n\geq0}$ is strongly $q$-log-convex, then so is
$y_n(q)=\sum_{k\geq0} T_{n,k}x_k(q)$. In particular, if
$\{x_n\}_{n\geq0}$ is log-convex, then so is $y_n=\sum_{k\geq0}
T_{n,k}x_k$.
\end{thm}

 For
$a_2\geq0$ and $b_2\geq0$, Chen {\it et al.} \cite{CWY11} proved
that $[T_{n,k}]_{n,k\ge 0}$ satisfying the recurrence
(\ref{rr+Chen}) is TP$_2$. For $a_2\geq0$ and $b_2\geq0$, we
\cite{Zhu14} also showed that $[T_{n,k}]_{n,k\ge 0}$ satisfying the
recurrence (\ref{rr+Zhu}) is TP. Thus the following result is
immediate from Theorem \ref{thm+transf+squar}.

\begin{thm}\label{thm+transf+chen+zhu}
Let $[T_{n,k}]_{n,k\geq0}$ be an array of nonnegative numbers
satisfying the recurrence relation:
\begin{equation}\label{rr+squar}
T_{n,k}=(a_0n+a_1k^2+a_2k+a_3)T_{n-1,k}+(b_0n+b_1k^2+b_2k+b_3)T_{n-1,k-1}
\end{equation}
with $T_{n,k}=0$ unless $0\le k\le n$ and $T_{0,0}=1$. Assume that
$a_2\geq0,b_2\geq0$ and $a_1=b_1=0$, or $a_2\geq0,b_2\geq0$ and
$a_0=b_0=0$. If $\{x_n(q)\}_{n\geq0}$ is strongly $q$-log-convex,
then so is $y_n(q)=\sum_{k=0}^n T_{n,k}x_k(q)$. In particular, if
$\{x_n\}_{n\geq0}$ is log-convex, then so is $y_n=\sum_{k=0}^n
T_{n,k}x_k$.
\end{thm}

\begin{rem}
For $a_2\geq0,b_2\geq0$ and $a_1=b_1=0$ in Theorem
\ref{thm+transf+chen+zhu}, recently Liu and Li \cite{LL16}
independently proved the result.
\end{rem}

 Twenty
five years ago Gian-Carlo Rota said ``The one contribution of mine
that I hope will be remembered has consisted in just pointing out
that all sorts of problems of combinatorics can be viewed as
problems of the locations of zeros of certain polynomials...'', see
the end of the introduction of \cite{BBL09}. In fact, polynomials
with only real zeros play an important role in attacking
log-concavity of sequences. One classical result is that if the
polynomial $\sum_{i=0}^{n}a_ix^i$ with nonnegative coefficients has
only real zeros, then the sequence $a_0,a_1,\ldots,a_n$ is
log-concave. In addition, many log-concave sequences arising in
combinatorics have the stronger property, see Liu and
Wang~\cite{LW-RZP} and Wang and Yeh~\cite{WYjcta05} for instance.
Using the algebraical method, Liu and Wang~\cite{LW07} found that
many polynomials with real zeros have the $q$-log-convexity. Thus at
the end of their paper, they proposed the problem to research this
relation between the $q$-log-convexity and real zeros. This is our
another motivation.

One of the classical problems of the theory of equations is to find
relations between the zeros and coefficients of a polynomial. A real
polynomial is {\it (Hurwitz) stable} if all of its zeros lie in the
open left half of the complex plane. A well-known necessary
condition for a real polynomial with positive leading coefficient to
be stable is that all its coefficients are positive. Polynomial
stability problems of various types arise in a number of problems in
mathematics and engineering. We refer to \cite[Chapter 9]{M66} for
deep surveys on the stability theory. Clearly, the stability
property of iterated polynomials implies the $q$-log-convexity. Thus
it is natural to consider the following stronger problem.
\begin{pro}\label{pro} Given a sequence $\{f_n(q)\}_{n\geq 0}$
of polynomials with only real zeros, under which conditions can we
obtain that $ f_{n+1}(q)f_{n-1}(q)-f^2_n(q)$ is stable for $n\geq1$
?
\end{pro}
We say a polynomial is a {\it generalized stable} polynomial if all
of its zeros excluding $0$ lie in the open left half of the complex
plane. The following result gives an answer to Problem \ref{pro}.
\begin{thm}\label{thm+poly+stable}
Let $\{f_n(q)\}_{n\geq0}$ be a sequence of  polynomials with
nonnegative coefficients, where $deg(f_n(q))=deg(f_{n-1}(q))+1$ for
$n\geq1$. Assume that the sequence $\{f_{n}(q)\}_{n\geq0}$ satisfies
the recurrence relation
$$f_n(q)=[a_1n+a_2+(b_1n+b_2)q+(c_1n+c_2)q^2]f_{n-1}(q)+q(a_3+b_3q+c_3q^2)f'_{n-1}(q),$$
where $a_1,b_1,c_1,a_1+a_3,b_1+b_3,c_1+c_3$ are all nonnegative. If
$\{f_{n}(q)\}_{n\geq 0}$ is a generalized Sturm sequence, then it is
$q$-log-convex. Furthermore, assume that  $c_1=c_3=0$. If
$a_1+2a_3\geq0$ and $b_1\geq b_3$, then $
f_{n+1}(q)f_{n-1}(q)-f^2_n(q)$ is a generalized stable polynomial
for $n\geq1$.
\end{thm}
The generalized Sturm sequences arise often in combinatorics. In
addition, the following result given by Liu and Wang~\cite{LW-RZP}
provides an approach to the generalized Sturm sequences.
\begin{prop}\label{prop+sturm}
Let $\{P_n(x)\}$ be a sequence of polynomials with nonnegative
coefficients and $\deg P_n=\deg P_{n-1}+1$. Suppose that
$$P_n(x)=(a_nx+b_n)P_{n-1}(x)+x(c_nx+d_n)P'_{n-1}(x)$$
where $a_n,b_n\in\mathbb{R}$ and $c_n\le 0,d_n\ge 0$. Then
$\{P_n(x)\}_{n\geq0}$ is a generalized Sturm sequence.
\end{prop}
It is well-known that many classical combinatorial sequences of
polynomials arising in certain triangular
 arrays, {\it e.g.}, Pascal triangle, Stirling triangle,
 Eulerian triangle and so on, satisfy the recurrence relation
 (\ref{rr+Chen}).
For its row generating functions $T_n(q)=\sum_{k=0}^n T_{n,k}q^k$,
by the recurrence relation (\ref{rr+Chen}), we have
$$T_n(q)=\left[a_1n+a_3+(b_1n+b_2+b_3)q\right]T_{n-1}(q)+(a_2+b_2q)qT'_{n-1}(q).$$
By Proposition~\ref{prop+sturm}, we know that if $a_2\ge0\geq b_2$
then the polynomials $T_n(q)$ form a generalized Sturm sequence.
Thus, the next result follows from Theorem~\ref{thm+poly+stable}.

\begin{prop}\label{prop+stable}
Let $[T_{n,k}]_{n,k\geq0}$ be the nonnegative array defined in
(\ref{rr+Chen}) and the row generating functions
$T_n(q)=\sum_{k=0}^n T_{n,k}q^k$. If $a_2\geq 0\geq b_2$, then
$\{T_{n+1}(q)T_{n-1}(q)-T_n^2(q)\}_{n\geq 1}$ is a sequence of
generalized stable polynomials.
\end{prop}

 The remainder of this paper is structured as follows.  In Section $2$, we will present the proofs of Theorem \ref{thm+transf+squar}.
 In Section $3$, we give the proof of
Theorem \ref{thm+poly+stable}. In Section $4$, we apply Theorem
\ref{thm+transf+squar} to some famous triangular arrays in a unified
manner, including Stirling triangles of two kinds, the
Jacobi-Stirling triangles of two kinds, the Legendre-Stirling
triangles of two kinds, the central factorial numbers triangle, the
Ramanujan transformation, and so on. In particular, we solve the
Conjecture \ref{Jacobi}. Finally, we also apply Proposition
\ref{prop+stable} to Eulerian polynomials of type $A$ and $B$, and
their $q$-analogs. Using Theorem \ref{thm+poly+stable}, we also
obtain the $q$-log-convexity of the generating functions of
alternating runs, the longest alternating subsequence and up-down
runs of permutations, respectively. In the Section $5$, we give some
remarks about linear transformations preserving the strong
$q$-log-convexity. In addition, we also present a criterion for the
strong $q$-log-convexity.

\section{Proof of Theorem \ref{thm+transf+squar}}


The next lemma plays an important role in our proof.
\begin{lem}\label{lem+equat}\cite{W02}
Given four sequences $\{a_i\}_{i=0}^n$, $\{b_i\}_{i=0}^n$,
$\{c_i\}_{i=0}^n$ and $\{d_i\}_{i=0}^n$, we have
$$\sum_{i=0}^na_ic_i\sum_{i=0}^nb_id_i-\sum_{i=0}^na_id_i\sum_{i=0}^nb_ic_i=\sum_{0\leq i<j\leq n}(a_ib_j-a_jb_i)(c_id_j-c_jd_i).$$
\end{lem}


 \textbf{Proof of Theorem
\ref{thm+transf+squar}:} In the following proof, we simply write
$x_k$ for $x_k(q)$.

In order to prove the strong $q$-log-convexity of
$\{y_n(q)\}_{n\geq0}$, it suffices to show for $n\geq m\geq1$ that
$$y_{n+1}(q)y_{m-1}(q)-y_{n}(q)y_{m}(q)\geq_q 0.$$
Then, for $n\geq m\geq1$, by the recurrence relation (\ref{rr+Gen}),
we have
\begin{eqnarray}
&&y_{n+1}(q)y_{m-1}(q)-y_{n}(q)y_{m}(q)\nonumber\\
&=&\sum_{k\geq0}[r(n+1)+f(k)]T_{n,k}x_k\sum_{k\geq0}T_{m-1,k}x_k-\sum_{k\geq0}T_{n,k}x_k\sum_{k\geq0}[r(m)+f(k)]T_{m-1,k}x_k+\nonumber\\
&&\sum_{k\geq1}[s(n+1)+g(k)]T_{n,k-1}x_k\sum_{k\geq0}T_{m-1,k}x_k-\sum_{k\geq0}T_{n,k}x_k\sum_{k\geq1}[s(m)+g(k)]T_{m-1,k-1}x_k+\nonumber\\
&&\sum_{k\geq2}[t(n+1)+h(k)]T_{n,k-2}x_k\sum_{k\geq0}T_{m-1,k}x_k-\sum_{k\geq0}T_{n,k}x_k\sum_{k\geq2}[t(m)+h(k)]T_{m-1,k-2}x_k\nonumber\\
&=&\left[\sum_{k\geq0}r(n+1)T_{n,k}x_k\sum_{k\geq0}T_{m-1,k}x_k-\sum_{k\geq0}T_{n,k}x_k\sum_{k\geq0}r(m)T_{m-1,k}x_k\right]+\label{sum1}\\
&&\left[\sum_{k\geq0}f(k)T_{n,k}x_k\sum_{k\geq0}T_{m-1,k}x_k-\sum_{k\geq0}T_{n,k}x_k\sum_{k\geq0}f(k)T_{m-1,k}x_k\right]+\label{sum2}\\
&&\left[\sum_{k\geq1}s(n+1)T_{n,k-1}x_k\sum_{k\geq0}T_{m-1,k}x_k-\sum_{k\geq0}T_{n,k}x_k\sum_{k\geq1}s(m)T_{m-1,k-1}x_k\right]+\label{sum3}\\
&&\left[\sum_{k\geq1}g(k)T_{n,k-1}x_k\sum_{k\geq0}T_{m-1,k}x_k-\sum_{k\geq0}T_{n,k}x_k\sum_{k\geq1}g(k)T_{m-1,k-1}x_k\right]+\label{sum4}\\
&&\left[\sum_{k\geq2}t(n+1)T_{n,k-2}x_k\sum_{k\geq0}T_{m-1,k}x_k-\sum_{k\geq0}T_{n,k}x_k\sum_{k=0}^mt(m)T_{m-1,k-2}x_k\right]+\label{sum5}\\
&&\left[\sum_{k\geq2}h(k)T_{n,k-2}x_k\sum_{k\geq0}T_{m-1,k}x_k-\sum_{k\geq0}T_{n,k}x_k\sum_{k=0}^mh(k)T_{m-1,k-2}x_k\right].\label{sum6}
\end{eqnarray}

In what follows we will prove that every difference in
(\ref{sum1})-(\ref{sum6}) is $q$-nonnegative.
 Obviously, for (\ref{sum1}), we have
\begin{eqnarray*}
&&\sum_{k\geq0}r(n+1)T_{n,k}x_k\sum_{k\geq0}T_{m-1,k}x_k-\sum_{k\geq0}T_{n,k}x_k\sum_{k\geq0}r(m)T_{m-1,k}x_k\nonumber\\
&=&[r(n+1)-r(m)]\sum_{k\geq0}T_{n,k}x_k\sum_{k\geq0}T_{m-1,k}x_k\nonumber\\
&\geq_q&0
\end{eqnarray*}
since $r(n)$ is nonnegative and increasing in $n$.

For (\ref{sum2}), if we view $T_{m-1,k}$, $T_{n,k}$, $x_k$ and
$f(k)x_k$ as $a_k$, $b_k$, $c_k$ and $d_k$ in Lemma \ref{lem+equat},
respectively, then we obtain that
\begin{eqnarray*}
&&\sum_{k\geq0}f(k)T_{n,k}x_k\sum_{k\geq0}T_{m-1,k}x_k-\sum_{k\geq0}T_{n,k}x_k\sum_{k\geq0}f(k)T_{m-1,k}x_k\nonumber\\
&=&\sum_{k=0}^{p}f(k)T_{n,k}x_k\sum_{k=0}^{p}T_{m-1,k}x_k-\sum_{k=0}^{p}T_{n,k}x_k\sum_{k=0}^{p}f(k)T_{m-1,k}x_k\nonumber\\
&=&\sum_{0\leq i<j\leq
p}[f(j)-f(i)](T_{m-1,i}T_{n,j}-T_{m-1,j}T_{n,i})x_ix_j\nonumber\\
&\geq_q&0
\end{eqnarray*}
since $f(k)$ is increasing and the matrix $[T_{n,k}]_{n,k\geq0}$ is
TP$_2$, where $$p=\max\{k:T_{n,k}\neq0\, or\, T_{m-1,k}\neq0\}.$$

Similarly by Lemma \ref{lem+equat}, for (\ref{sum3}) and
(\ref{sum4}), we derive that
\begin{eqnarray*}
&&\sum_{k\geq1}T_{n,k-1}x_k\sum_{k\geq0}T_{m-1,k}x_k-\sum_{k\geq0}T_{n,k}x_k\sum_{k\geq1}T_{m-1,k-1}x_k\nonumber\\
&=&\sum_{k\geq0}T_{n,k}x_{k+1}\sum_{k\geq0}T_{m-1,k}x_k-\sum_{k\geq0}T_{n,k}x_k\sum_{k\geq0}T_{m-1,k}x_{k+1}\nonumber\\
&=&\sum_{k=0}^pT_{n,k}x_{k+1}\sum_{k=0}^pT_{m-1,k}x_k-\sum_{k=0}^pT_{n,k}x_k\sum_{k=0}^pT_{m-1,k}x_{k+1}\nonumber\\
&=&\sum_{0\leq i<j\leq
p}\left[T_{m-1,i}T_{n,j}-T_{m-1,j}T_{n,i}\right]\left[x_ix_{j+1}-x_{i+1}x_j\right]\nonumber\\
&\geq_q&0
\end{eqnarray*}
and
\begin{eqnarray*}
&&\sum_{k\geq1}g(k)T_{n,k-1}x_k\sum_{k\geq0}T_{m-1,k}x_k-\sum_{k\geq0}T_{n,k}x_k\sum_{k\geq1}g(k)T_{m-1,k-1}x_k\nonumber\\
&=&\sum_{k\geq0}g(k+1)T_{n,k}x_{k+1}\sum_{k\geq0}T_{m-1,k}x_k-\sum_{k\geq0}T_{n,k}x_k\sum_{k\geq0}g(k+1)T_{m-1,k}x_{k+1}\nonumber\\
&=&\sum_{k=0}^pg(k+1)T_{n,k}x_{k+1}\sum_{k=0}^pT_{m-1,k}x_k-\sum_{k=0}^pT_{n,k}x_k\sum_{k=0}^pg(k+1)T_{m-1,k}x_{k+1}\nonumber\\
&=&\sum_{0\leq i<j\leq
p}\left[T_{m-1,i}T_{n,j}-T_{m-1,j}T_{n,i}\right]\left[g(j+1)x_ix_{j+1}-g(i+1)x_{i+1}x_j\right]\nonumber
\end{eqnarray*}
\begin{eqnarray*}
&\geq&\sum_{0\leq i<j\leq
n}\left[T_{m-1,i}T_{n,j}-T_{m-1,j}T_{n,i}\right]g(j+1)\left[x_ix_{j+1}-x_{i+1}x_j\right]\nonumber\\
&\geq_q&0.
\end{eqnarray*}
For (\ref{sum5}) and (\ref{sum6}), we obtain that
\begin{eqnarray*}
&&\sum_{k\geq2}T_{n,k-2}x_k\sum_{k\geq0}T_{m-1,k}x_k-\sum_{k\geq0}T_{n,k}x_k\sum_{k\geq2}T_{m-1,k-2}x_k\nonumber\\
&=&\sum_{k\geq0}T_{n,k}x_{k+2}\sum_{k\geq0}T_{m-1,k}x_k-\sum_{k\geq0}T_{n,k}x_k\sum_{k\geq0}T_{m-1,k}x_{k+2}\nonumber\\
&=&\sum_{k=0}^pT_{n,k}x_{k+2}\sum_{k=0}^pT_{m-1,k}x_k-\sum_{k=0}^pT_{n,k}x_k\sum_{k=0}^pT_{m-1,k}x_{k+2}\nonumber\\
&=&\sum_{0\leq i<j\leq
p}\left[T_{m-1,i}T_{n,j}-T_{m-1,j}T_{n,i}\right]\left[x_ix_{j+2}-x_{i+2}x_j\right]\nonumber\\
&\geq_q&0
\end{eqnarray*}
and
\begin{eqnarray*}
&&\sum_{k\geq2}h(k)T_{n,k-2}x_k\sum_{k\geq0}T_{m-1,k}x_k-\sum_{k\geq0}T_{n,k}x_k\sum_{k\geq2}h(k)T_{m-1,k-2}x_k\nonumber\\
&=&\sum_{k\geq0}h(k+2)T_{n,k}x_{k+2}\sum_{k\geq0}T_{m-1,k}x_k-\sum_{k\geq0}T_{n,k}x_k\sum_{k\geq0}h(k+2)T_{m-1,k}x_{k+2}\nonumber\\
&=&\sum_{k=0}^ph(k+2)T_{n,k}x_{k+2}\sum_{k=0}^pT_{m-1,k}x_k-\sum_{k=0}^pT_{n,k}x_k\sum_{k=0}^ph(k+2)T_{m-1,k}x_{k+2}\nonumber\\
&=&\sum_{0\leq i<j\leq
p}\left[T_{m-1,i}T_{n,j}-T_{m-1,j}T_{n,i}\right]\left[h(j+2)x_ix_{j+2}-h(i+2)x_{i+2}x_j\right]\nonumber\\
&\geq&\sum_{0\leq i<j\leq
p}\left[T_{m-1,i}T_{n,j}-T_{m-1,j}T_{n,i}\right]h(i+2)\left[x_ix_{j+2}-x_{i+2}x_j\right]\nonumber\\
&\geq_q&0.
\end{eqnarray*}
 This completes the proof.
\qed

\section{Proof of Theorem \ref{thm+poly+stable}}
\hspace*{\parindent} Following Wagner~\cite{Wag92}, a real
polynomial is said to be {\it standard} if either it is identically
zero or its leading coefficient is positive. Suppose that both $f$
and $g$ only have real zeros. Let $\{r_i\}$ and $\{s_j\}$ be all
zeros of $f$ and $g$ in nondecreasing order respectively. We say
that $g$ {\it interlaces} $f$ if $\deg f=\deg g+1=n$ and
\begin{equation}\label{int-def}
r_n\le s_{n-1}\le\cdots\le s_2\le r_2\le s_1\le r_1.
\end{equation}
By $g\sep f$ we denote ``$g$ interlaces $f$''. For notational
convenience, let $a\sep bx+c$ for any real constants $a,b,c$ and
$f\sep 0, 0\sep f$ for all real polynomial $f$ with only real zeros.

Let $\{P_n(x)\}_{n\ge 0}$ be a sequence of standard polynomials.
Recall that $\{P_n(x)\}$ is a {\it Sturm sequence} if $\deg P_n=n$
and $P_{n-1}(r)P_{n+1}(r)<0$ whenever $P_n(r)=0$ and $n\ge 1$. We
say that $\{P_n(x)\}$ is a {\it generalized Sturm sequence} if
$P_n\in\rz$ and $P_0\sep P_1\sep\cdots\sep P_{n-1}\sep
P_n\sep\cdots$. For example, if $P$ is a standard polynomial with
only real zeros and $\deg P=n$, then $P^{(n)},P^{(n-1)},\ldots,P',P$
form a generalized Sturm sequence by Rolle's theorem.

In order to simplify our proof, we need the following lemma.
\begin{lem}\emph{\cite[Lemma 1.20]{Fisk08}}\label{lem+Fisk}
Let both $f(x)$ and $g(x)$ be standard real polynomials with only
real zeros. Assume that $\deg(f(x))=n$ and all real zeros of $f(x)$
are $s_1, \ldots, s_n$. If $\deg(g)=n-1$ and we write
$$g(x)=\sum_{i=1}^n\frac{c_if(x)}{x-s_i},$$
then $g(x)$ interlaces $f(x)$ if and only if all $c_i$ are positive.
\end{lem}

\textbf{Proof of Theorem \ref{thm+poly+stable}:}

Since
$$f_n(q)=[a_1n+a_2+(b_1n+b_2)q+(c_1n+c_2)q^2]f_{n-1}(q)+q(a_3+b_3q+c_3q^2)f'_{n-1}(q),$$
it follows that
\begin{eqnarray}
&&f_{n+1}(q)f_{n-1}(q)-f_{n}^2(q)\nonumber\\
&=&[a_1n+a_2+a_1+(b_1n+b_2+b_1)q+(c_1n+c_2+c_1)q^2]f_{n}(q)f_{n-1}(q)+\nonumber\\
&&q(a_3+b_3q+c_3q^2)f'_{n}(q)f_{n-1}(q)-q(a_3+b_3q+c_3q^2)f'_{n-1}(q)f_{n}(q)-\nonumber\\
&&[a_1n+a_2+(b_1n+b_2)q+(c_1n+c_2)q^2]f_{n-1}(q)f_{n}(q) \nonumber\\
&=&(a_1+b_1q+c_1q^2)f_{n}(q)f_{n-1}(q)+q(a_3+b_3q+c_3q^2)\left[f'_{n}(q)f_{n-1}(q)-f_{n}(q)f'_{n-1}(q)\right]\nonumber\\
&=&f^2_n(q)\left[(a_1+b_1q+c_1q^2)\frac{f_{n-1}(q)}{f_{n}(q)}-q(a_3+b_3q+c_3q^2)\left(\frac{f_{n-1}(q)}{f_{n}(q)}\right)'\right]\label{poly}.
\end{eqnarray}
By assumption, $\{f_{n}(q)\}_{n\geq 0}$ is a generalized Sturm
sequence. Thus, if we assume that the all non-positive zeros of
$f_n(q)$ are ordered as $r_1\ge r_2 \ge\ldots \ge r_n$, then
$f_{n-1}(q)=f_n(q)\sum_{i=1}^n\frac{s_i}{q-r_i}$ by
Lemma~\ref{lem+Fisk}, where $s_i> 0$ for $1\leq i\leq n$. Hence,
\begin{eqnarray*}\label{qq}
&&f_{n+1}(q)f_{n-1}(q)-f_{n}^2(q)\\
&=&f^2_n(q)\left[(a_1+b_1q+c_1q^2)\frac{f_{n-1}(q)}{f_{n}(q)}-q(a_3+b_3q+c_3q^2)\left(\frac{f_{n-1}(q)}{f_{n}(q)}\right)'\right]\\
&=&f^2_n(q)\left[(a_1+b_1q+c_1q^2)\sum_{i=1}^n\frac{s_i}{q-r_i}+q(a_3+b_3q+c_3q^2)\sum_{i=1}^n\frac{s_i}{(q-r_i)^2}\right]\\
&=&f^2_n(q)\sum_{i=1}^n\frac{s_i\left[(a_1+b_1q+c_1q^2)(q-r_i)+q(a_3+b_3q+c_3q^2)\right]}{(q-r_i)^2}\\
&=&\sum_{i=1}^ns_i\left[(c_1+c_3)q^3+(b_1+b_3-c_1r_i)q^2+(a_1+a_3-b_1r_i)q-a_1r_i\right]\left(\frac{f_n(q)}{q-r_i}\right)^2,
\end{eqnarray*}
which is a polynomial with nonnegative coefficients since
$$(c_1+c_3)q^3+(b_1+b_3-c_1r_i)q^2+(a_1+a_3-b_1r_i)q-a_1r_i, \frac{f_n(q)}{q-r_i}$$
are all polynomials with nonnegative coefficients for $1\leq i \leq
n$. Thus $\{f_{n}(q)\}_{n\geq 0}$ is $q$-log-convex.

In the following, we proceed to demonstrate the second part that
$$f_{n+1}(q)f_{n-1}(q)-f_{n}^2(q)$$
is a generalized stable polynomial for each $n\geq 1$. Note that
\begin{eqnarray*}
&&f_{n+1}(q)f_{n-1}(q)-f_{n}^2(q)\\
&=&f^2_n(q)\sum_{i=1}^n\frac{s_i\left[(a_1+b_1q+c_1q^2)(q-r_i)+q(a_3+b_3q+c_3q^2)\right]}{(q-r_i)^2}.
\end{eqnarray*}
Thus we only need to show that
\begin{eqnarray}\label{sum}
\sum_{i=1}^n\frac{s_i\left[(b_1+b_3)q^2+(a_1+a_3-b_1r_i)q-a_1r_i\right]}{(q-r_i)^2}
\end{eqnarray}
has no zeros in the right half plane since $c_1=c_3=0$. Let
$q=x+yI$, where $I$ is the imaginary number unit. Then, for $x\geq0$
and $r\leq0$, it follows from $a_1+2a_3\geq0$ and $b_1\geq b_3$ that
we have
\begin{eqnarray*}
&&[(x-r)^2+y^2]^2Re\left(\frac{(b_1+b_3)q^2+(a_1+a_3-b_1r)q-a_1r}{(q-r)^2}\right)\\
&=&Re\left([(b_1+b_3)(x+yI)^2+(a_1+a_3-b_1r)(x+yI)-a_1r](x-r-yI)^2\right)\\
&=&B(x^2-y^2)^2+Bx^2(r^2-2xr)+(xA-xb_1r-a_1r)(x-r)^2+\\
&&y^2\left\{4x^2B+xA-(2B+b_1)xr-(a_1+2a_3)r+(b_1-b_3)r^2\right\}\\
&\geq&0,
\end{eqnarray*}
where $A=a_1+a_3$ and $B=b_1+b_3$. Thus,
$$Re\left(\sum_{k=1}^n\frac{c_k\left[(b_1+b_2)q^2+(a_1+a_2-b_1r_k)q-a_1r_k\right]}{(q-r_k)^2}\right)\geq0$$
with the equality if and only if $q=0$. This completes the proof.
\qed


\section{Applications}
In this section, we give some applications of the main results.
\subsection{Stirling transformations of two kinds}
Let $S_{n,k}$ denote the number of partitions of a set with $n$
elements consisting of $k$ disjoint nonempty sets. It is well known
that $S_{n,k}$ is called the Stirling number of the second kind. In
addition, the Stirling numbers of the second kind satisfy the
recurrence
\begin{equation*}
S_{n+1,k}=kS_{n,k}+S_{n,k-1}.
\end{equation*}
Let $c_{n,k}$ be the signless Stirling number of the first kind,
i.e., the number of permutations of $[n]$ which contain exactly $k$
permutation cycles. Similarly, signless Stirling numbers of the
first kind $c_{n,k}$ satisfy the recurrence
$$
c_{n,k}=(n-1)c_{n-1,k}+c_{n-1,k-1}.
$$

The $q$-log-convexity and strong $q$-log-convexity of the row
generating functions $S_n(q)=\sum_{k=0}^{n}S_{n,k}q^k$ have been
proved, see Liu and Wang~\cite{LW07}, Chen {\it et al.}~\cite{CWY11}
and Zhu \cite{Zhu13,Zhu14} for instance. Liu and Wang~\cite{LW07}
also proved that both the Stirling transformations of two kinds
preserve the log-convexity. By Theorem \ref{thm+transf+chen+zhu}, we
can extend above results to the strong $q$-log-convexity as follows.
\begin{prop} The linear transformation
$y_n(q)=\sum_{k=0}^{n}S_{n,k}x_k(q)$ preserves the strong
$q$-log-convexity.
\end{prop}
\begin{prop} The linear transformation
$y_n(q)=\sum_{k=0}^{n}c_{n,k}x_k(q)$ preserves the strong
$q$-log-convexity.
\end{prop}

\subsection{Jacobi-Stirling transformation of the second kind}

The Jacobi-Stirling numbers $\JS_n^k( z)$ of the second kind satisfy
the following recurrence relation:
\begin{equation}
\left\{ \begin{array}{l} \JS_0^0(z)=1, \qquad \JS_n^k(z)=0, \quad  \text{if} \ k \not\in\{1,\ldots,n\}, \\
\JS_n^k(z)= \JS_{n-1}^{k-1}(z)+k(k+z)\,\JS_{n-1}^{k}(z),  \quad n,k
\geq 1,\end{array} \right.
 \end{equation}
 where $z=\alpha+\beta+1.$ Lin and Zeng \cite{LZ14} and Zhu \cite{Zhu14} independently proved
the strong $q$-log-convexity of the row generating functions of
Jacobi-Stirling numbers. By Theorem \ref{thm+transf+chen+zhu}, we
have the following generalization, which in particular confirms
Conjecture \ref{Jacobi} for Jacobi-Stirling numbers of the second
kind.
\begin{prop}
The Jacobi-Stirling transformation of the second kind
$$y_n(q)=\sum_{k=0}^{n}\JS_n^k(z)x_k(q)$$ preserves the strong $q$-log-convexity for
$z\geq0$. In particular,
$$y_n=\sum_{k=0}^{n}\JS_n^k(z)x_k$$ preserves the log-convexity for
$z\geq0$.
\end{prop}

The Jacobi-Stirling numbers $\Jc_n^k( z)$ of the first kind are
defined by
\begin{equation*}
\left\{ \begin{array}{l} \Jc_0^0(z)=1, \qquad \Jc_n^k(z)=0, \quad  \text{if} \ k \not\in\{1,\ldots,n\}, \\
\Jc_n^k(z)= \Jc_{n-1}^{k-1}(z)+(n-1)(n-1+z)\,\Jc_{n-1}^{k}(z), \quad
n,k \geq 1,\end{array} \right.
 \end{equation*}
 where $z=\alpha+\beta+1.$
It is well known that the unsigned inverse of a totally positive
matrix is also totally positive. Thus by the total positivity of
$[\JS_n^k(z)]_{n,k}$ \cite{Zhu14}, we derive that
$[\Jc_n^k(z)]_{n,k}$ is totally positive, also see \cite{Monge12}.
Thus by Theorem \ref{thm+transf+squar}, we also have the following
result, which in particular confirms Conjecture \ref{Jacobi} for
Jacobi-Stirling numbers of the first kind.
\begin{prop}
The Jacobi-Stirling transformation of the first kind
$$w_n(q)=\sum_{k=0}^{n}\Jc_n^k(z)x_k(q)$$ preserves the strong $q$-log-convexity for
$z\geq0$. In particular,
$$w_n=\sum_{k=0}^{n}\JS_n^k(z)x_k$$ preserves the log-convexity for
$z\geq0$.
\end{prop}
\begin{rem}
If $z=1$, then $\JS_{n}^{k}(1)$ and $\Jc_{n}^{k}(1)$ are the
Legendre-Stirling numbers of two kinds, respectively.
\end{rem}

\subsection{Central factorial transformations}
The \emph{central factorial numbers} of the second kind $T(n,k)$ are
defined in Riordan's book \cite[p. 213-217]{Riordan} by
\begin{align}
x^n=\sum_{k=0}^nT(n,k)\,x\prod_{i=1}^{k-1}\left(x+\frac{k}{2}-i\right).
\end{align}
Therefore, if let $U(n,k)=T(2n,2k)$ and $ V(n,k)=4^{n-k}
T(2n+1,2k+1)$, then
\begin{eqnarray*}
U(n,k)&=&U(n-1,k-1)+k^2U(n-1,k),\\
V(n,k)& =& V(n-1,k-1) + (2k+1)^2 V(n-1,k).
 \end{eqnarray*}
 Zhu \cite{Zhu14} proved that the row generating functions of
 $U(n,k)$ (respectively, $V(n,k)$) form a strongly $q$-log-convex sequence.
  In view of Theorem
 \ref{thm+transf+chen+zhu}, these can be extended to the following result.
 \begin{prop} The linear transformation
$y_n(q)=\sum_{k=0}^{n}U(n,k)x_k(q)$ preserves the strong
$q$-log-convexity. In particular, $y_n=\sum_{k=0}^{n}U(n,k)x_k$
preserves the log-convexity.
\end{prop}
\begin{prop} The linear transformation
$y_n(q)=\sum_{k=0}^{n}V(n,k)x_k(q)$ preserves the strong
$q$-log-convexity. In particular, $y_n=\sum_{k=0}^{n}V(n,k)x_k$
preserves the log-convexity.
\end{prop}

\subsection{Ramanujan transformation}
Let $r_{n,k}$ be the number of rooted labeled trees on $n$ vertices
with $k$ improper edges. Then numbers $r_{n,k}$ satisfy the
following recurrence relation:
\begin{equation*}
r_{n,k}=(n-1)r_{n-1,k}+(n+k-2)r_{n-1, k-1}
\end{equation*}
where $r_{1,0}=1$, $n\geq1$, $k\leq n-1$, and $r_{n,k}=0$ otherwise,
see Shor \cite{Sh95}. It was proved that the row generating
functions of $[r_{n,k}]_{n,k\geq0}$ are the famous Ramanujan
polynomials $r_n(y)$, which are defined by the recurrence relation
$$r_1(y)=1, \ r_{n+1}=n(1+y) r_n(y)+y^2r_n'(y).$$
The first values of the polynomials $r_n(y)$ are
$$r_2(y)=1+y,\ r_3(y)=2+4y+3y^2,\ r_4(y)=6 +18y +25y^2+15y^3.$$
Chen {\it et al.} \cite{CWY11} proved that the polynomials $r_n(y)$
form a strongly $q$-log-convex sequence, which can be extended to
the following result by Theorem \ref{thm+transf+chen+zhu}.

\begin{prop} If $\{x_n(q)\}_{n\geq0}$ is strongly $q$-log-convex, then so is $\{y_n(q)\}_{n\geq0}$ defined by
$y_n(q)=\sum_{k=0}^{n}r_{n,k}x_k(q)$. In particular,
$y_n=\sum_{k=0}^{n}r_{n,k}x_k$ preserves the log-convexity.
\end{prop}

\subsection{Associated Lah transformation }
The associated Lah numbers defined by
$$L_m(n,k)=(n!/k!)\sum_{i=1}^k(-1)^{k-i}\binom{k}{i}\binom{n+mi-1}{n}$$
satisfy the recurrence
$$L_m(n,k)=(mk+n-1)L_m(n-1,k)+mL_m(n-1,k-1).$$
Let $L_n(q)=\sum_{k=0}^nL_m(n,k)q^k$. By virtue of Theorem
\ref{thm+transf+chen+zhu} and Proposition \ref{prop+stable},
 we have the following two results, respectively.
\begin{prop} The linear transformation
$z_n(q)=\sum_{k=0}^{n}L_m(n,k)x_k(q)$ preserves the strong
$q$-log-convexity. In particular, $z_n=\sum_{k=0}^{n}L_m(n,k)x_k$
preserves the log-convexity.
\end{prop}
\begin{prop}
 $L_{n+1}(q)L_{n-1}(q)-L_{n}^2(q)$ is a stable polynomial
 for each $n\geq 1.$
\end{prop}

\subsection{Eulerian polynomials of types $A$ and $B$}

Let $\pi=a_1a_2\cdots a_n$ be a permutation of $[n]$. An element
$i\in [n-1]$ is called a descent of $\pi$ if $a_i>a_{i+1}$. The
number of permutations of $[n]$ having $k-1$ descents is called the
Eulerian number, denoted by $A_{n,k}$, and its row-generating
function $A_n(q)=\sum_{k=0}^{n}A_{n,k}q^k$ is called the classical
Eulerian polynomial. It is known that
$$A_n(q)=nqA_{n-1}(q)+q(1-q)A'_{n-1}(q).$$
We refer reader to Comtet \cite{Com74} for furhter properties about
Eulerian polynomials. Let $B_{n,k}$ be the Eulerian number of type
$B$ counting the elements of $B_n$ with $k$ $B$-descents. It is
known that the Eulerian numbers of type $B$ satisfy the recurrence
\begin{equation}\label{recurrenc-1}B_{n,k}=(2k+1)B_{n-1,k}+(2n-2k+1)B_{n-1,k-1}.\end{equation}
Assume that $B_n(q)=\sum_{k=0}^{n}B_{n,k}q^k$ is the Eulerian
polynomial of type $B$. Then we have
$$B_n(q)=(1+q)B_{n-1}(q)+2x(1-x)B'_{n-1}(q).$$ It was proved that
polynomials $A_n(q)$ (respectively $B_n(q)$) form a strongly
$q$-log-convex sequence, respectively, see \cite{Zhu14,LZ15}. By
Theorem \ref{thm+poly+stable}, the following result is immediate.

\begin{prop} Both $A_{n+1}(q)A_{n-1}(q)-A^2_{n}(q)$ and $B_{n+1}(
q)B_{n-1}(q)-B^2_{n}(q)$ are generalized stable polynomials in $q$
for $n\geq1$.
\end{prop}

\begin{rem}
The polynomial $A_{n+1}(q)A_{n-1}(q)-A^2_n(q)$ is stable, which was
proved by Fisk \cite[Lemma 21.92]{Fisk08}.
\end{rem}

\subsection{q-Eulerian polynomials}

For a finite Coxeter group $W$, let $d_W(\pi)$ denote the number of
$W$-descents of $\pi$. Then the Eulerian polynomial of $W$ is
defined by
$$P(W,x)=\sum_{\pi\in W}x^{d_W(\pi)},$$
see Bj\"orner and Brenti \cite{BB05} for instance.

 For Coxeter groups of type $A$,
$P(A_n,x)=A_n(x)/x$, where $A_n(x)$ is the classical Eulerian
polynomial. Let $\exc(\pi)$ and $c(\pi)$ denote the numbers of
excedances and cycles in $\pi$, respectively. In \cite{FS70}, Foata
and Sch\"utzenberger defined a $q$-analog of the classical Eulerian
polynomials by
$$A_n(x;q)=\sum_{\pi\in\msn}x^{\exc(\pi)+1}q^{c(\pi)}.$$
Obviously, for $q=1$, $A_n(x;q)$ reduces to the classical Eulerian
polynomial $A_n(x)$. In addition, in \cite[Proposition 7.2]{Bre00},
Brenti demonstrated the recurrence relation
$$A_{n}(x;q)=(nx+q-1)A_{n-1}(x;q)+x(1-x)\frac{\partial}{\partial x}A_{n-1}(x;q),$$
with the initial condition $A_0(x;q)=x$.

For Coxeter groups of type $B$, letting $N(\pi)=|\{i\in [n]:
\pi(i)<0\}|$, Brenti~\cite{Bre94EuJC} introduced a $q$-analogue of
$P(B_n,x)$ by
$$B_n(x;q)=\sum_{\pi\in B_n}q^{N(\pi)}x^{d_B(\pi)}.$$
In fact, $B_n(x;q)$ reduces to $A_n(x)$ for $q=0$ and to $P(B_n,x)$
for $q=1$. It is also known that $\{B_n(x;q)\}_{n\geq0}$ satisfies
the recurrence relation
$$B_n(x;q)=\{1+[(1+q)n-1]x\}B_{n-1}(x;q)+(1+q)x(1-x)\frac{\partial}{\partial x}B_{n-1}(x;q),$$
with $B_0(x;q)=1$, see \cite[Theorem 3.4 (i)]{Bre94EuJC}. Thus the
following result follows from Proposition \ref{prop+stable}.
\begin{prop}
Both $A_{n+1}(x; q)A_{n-1}(x;q)-A^2_{n}(x;q)$ and $B_{n+1}(x;
q)B_{n-1}(x;q)-B^2_{n}(x;q)$ are generalized stable polynomials in
$x$ for any fixed $q\geq0$.
\end{prop}

\subsection{Alternating runs}
Suppose that $\msn$ denotes the symmetric group of all permutations
of $\{1,2,\ldots,n\}$. For $\pi=\pi(1)\pi(2)\cdots \pi(n)\in\msn$,
we say that $\pi$ changes direction at position $i$ if either
$\pi({i-1})<\pi(i)>\pi(i+1)$, or $\pi(i-1)>\pi(i)<\pi(i+1)$. We say
that $\pi$ has $k$ alternating runs if there are $k-1$ indices $i$
such that $\pi$ changes direction at these positions. Denote by
$R(n, k)$ the number of permutations in $S_n$ having $k$ alternating
runs. Then we have
\begin{eqnarray}\label{run}
R(n, k) = k\,R(n- 1, k)+2R(n-1, k -1) + (n- k)R(n-1, k-2)
\end{eqnarray}  for $n, k\geq1$, where $R(1, 0) = 1$ and $R(1, k) =
0$ for $k\geq1$, see B\'{o}na \cite{Bona12} for instance. For
$n\geq1$, define the alternating runs polynomials $R_n(x) =
\sum_{k=1}^{n-1}R(n, k)x^k$. It follows from the recurrence
(\ref{run}) that
\begin{eqnarray*}
R_{n+2}(x) = x(nx + 2)R_{n+1}(x) + x(1-x^2)R'_{n+1}(x)
\end{eqnarray*}
with initial conditions $R_1(x) = 1$ and $R_2(x) = 2x$. Moreover,
the polynomial $R_n(x)$ has a close connection with the classical
Eulerian polynomial $A_n(x)$ by
\begin{eqnarray*}
R_n(x) = \left(\frac{1+x}{2}\right)^{n-1}
(1+w)^{n+1}A_n(\frac{1-w}{1+w}), w =\sqrt{\frac{1 - x}{1 + x}},
\end{eqnarray*} see Knuth \cite{Kur73}. The polynomials $R_n(x)$ also have only non-positive real zeros and
$R_n(x)\sep R_{n+1}(x)$, see Ma and Wang~\cite{MW08}. It follows
from Theorem \ref{thm+poly+stable} that the next result is
immediate.
\begin{prop}
 The alternating runs polynomials $R_n(q)$ form a $q$-log-convex sequence.
\end{prop}
\subsection{The longest alternating subsequence and up-down runs of
permutations}

For a subsequence $\pi({i_1})\cdots \pi({i_k})$ of $\pi$, it is
called an {\it alternating subsequence} if
$$\pi({i_1})>\pi({i_2})<\pi({i_3})>\cdots \pi({i_k}).$$
Let $\as(\pi)$ and and $a_k(n)$ denote the length of the longest
alternating subsequence of $\pi$ and the number of permutations in
$\msn$ with $\as(\pi)=k$. Define its ordinary generating function
$t_n(x)=\sum_{k=1}^na_k(n)x^k$. For $n\geq 2$, B\'ona~\cite[Section
1.3.2]{Bona12} derived the following identity:
\begin{equation*}
t_n(x)=\frac{1}{2}(1+x)R_n(x).
\end{equation*}
Ma \cite{Ma13} also proved that the polynomials $t_n(x)$ satisfy the
recurrence relation
\begin{equation*}\label{Tnx-recurrence}
t_{n+1}(x)=x(nx+1)t_{n}(x)+x\left(1-x^2\right)t_{n}'(x),
\end{equation*}
with initial conditions $t_0(x)=1$ and $t_1(x)=x$. We refer reader
to Stanley \cite{Sta08} for more properties about the longest
alternating subsequences.

On the other hand, $a_k(n)$ is also the number of permutations in
$\msn$ with $k$ up-down runs. The {\it up-down runs} of a
permutation $\pi$ are defined to the alternating runs of $\pi$
endowed with a $0$ in the front, see~\cite[A186370]{Sloane}. In
addition, the up-down runs of a permutation have a close connection
with interior peaks and left peaks. Based on the interior peaks and
left peaks, Ma \cite{Ma121} defined polynomials $M_n(x)$, which
satisfy the recurrence relation
\begin{equation*}
M_{n+1}(x)=(1+nx^2)M_{n}(x)+x(1-x^2)M_n'(x),
\end{equation*}
with initial conditions $M_1(x)=1+x$ and $M_2(x)=1+2x+x^2$,
see~\cite[Section 2]{Ma121}. In addition, $M_n(x)$ has only
non-positive real zeros and $M_n(x)\sep M_{n+1}(x)$, see Ma
\cite{Ma121}.
So, we have the following result by Theorem \ref{thm+poly+stable}.
\begin{prop}
 Both $\{t_n(q)\}_{n\geq 0}$ and $\{M_n(q)\}_{n\geq 0}$ are $q$-log-convex sequences.
\end{prop}

\subsection{Alternating runs of type $B_n$}
A {\it run} of a signed permutation $\pi\in B_n$ is defined as a
maximal interval of consecutive elements on which the elements of
$\pi$ are monotonic in the order
$\cdots<\overline{2}<\overline{1}<0<1<2<\cdots$. Let $T(n,k)$ denote
the number of signed permutations in $B_n$ with $k$ alternating runs
and $\pi(1)> 0$. In \cite[Theorem 4.2.1]{Zhao11}, it was shown that
the array $[T(n,k)]_{n,k}$ satisfies the recurrence relation
\begin{equation}\label{tnk-recurrence03}
T(n,k)=(2k-1)T(n-1,k)+3T(n-1,k-1)+(2n-2k+2)T(n-1,k-2)
\end{equation}
for $n\geqslant 2$ and $1\leqslant k\leqslant n$, where $T(1,1)=1$
and $T(1,k)=0$ for $k>1$. Let $T_n(x)=\sum_{k=1}^nT(n,k)x^k$ denote
the alternating run polynomials of type $B$. It follows from
(\ref{tnk-recurrence03}) that we have the recurrence relation
\begin{align*}
T_n(x)=[2(n-1)x^2+3x-1]T_{n-1}(x)+2x(1-x^2)T'_{n-1}(x).
\end{align*}
Zhao proved that polynomials $T_n(x)$ form a generalized Sturm
sequence by using Theorem $2$ of Ma and Wang \cite{MW08}. Thus we
get the following result from Theorem \ref{thm+poly+stable}.
\begin{prop}
For $n\geq1$, the alternating run polynomials $T_n(q)$ form a
$q$-log-convex sequence.
\end{prop}

\section{Remarks}
There are also many famous triangular arrays, including the Motzkin
triangle, the Bell triangle, the Catalan triangle, the large
Schr\"{o}der triangle, and so on, satisfying such recurrence
relation
\begin{equation}\label{rec-three+right}
T_{n,k}=f_kT_{n-1,k-1}+g_kT_{n-1,k}+h_{k}T_{n-1,k+1}
\end{equation}
with $T_{0,0}=1$ and $T_{n,k}=0$ unless $0\leq k\leq n$ and the
nonnegative array $[T_{n,k}]_{n,k}$ has a general combinatorial
interpretation from the weighted Motzkin path,  see \cite{Fla80}. In
\cite{Zhu13}, we proved that the array $[T_{n,k}]_{n,k}$ is TP$_2$
and the first column $\{T_{n,0}\}_{n\geq0}$ is log-convex if
$g_{k+1}g_{k}\geq h_kf_{k+1}$ for $k\geq0$. For this array
$[T_{n,k}]_{n,k}$ in (\ref{rec-three+right}), by Theorem
\ref{thm+transf+squar}, it is natural to ask whether we have a
similar result. However, for the general case, the answer is not. In
the following, we give a simple example.

\begin{ex}\label{example}
Let a triangular array $[M_{n,k}]_{n,k}$ satisfy the recurrence
relation
$$M_{n,k}=M_{n-1,k-1}+M_{n-1,k}+M_{n-1,k+1}$$
with $M_{0,0}=1$ and $M_{n,k}=0$ unless $0\leq k\leq n$. This array
is called the Motzkin triangle and $M_{n,0}$ is the Motzkin number.
Assume that $x_0(q)=1$ and $x_n(q)=2^{n-1}q^n$ for $n\geq1$. It is
obvious that $\{x_n(q)\}_{n\geq0}$ is strongly $q$-log-convex. Let
$y_n=\sum_{k\geq0}M_{n,k}x_k(q)$ for $n\geq0$. It is easy to get
$$y_3y_1-y_2^2=q-q^2+2q^3,$$
which is not $q$-nonnegative. Thus the transformation
$y_n=\sum_{k\geq0}M_{n,k}x_k(q)$ does not preserve the strong
$q$-log-convexity.
\end{ex}
If all $f_k$, $g_k$ and $h_k$  in the (\ref{rec-three+right}) are
constants, then in \cite{Zhu17}, we gave a result for the strong
$q$-log-convexity of its row-generating functions. In fact, using
Lemma \ref{lem+equat}, similar to the proof of Theorem
\ref{thm+transf+squar}, we can get the following generalized
criterion, whose proof is omitted for brevity.

\begin{thm}
Let $\{f_n\}_{n\geq0}$, $\{g_n\}_{n\geq0}$ and $\{h_n\}_{n\geq0}$ be
nonnegative and increasing sequences, respectively. Define a
triangular array $[T_{n,k}]_{n,k\geq0}$ by
\begin{eqnarray*}\label{recc}
T_{n,k}=f_kT_{n-1,k-1}+g_kT_{n-1,k}+h_{k}T_{n-1,k+1}
\end{eqnarray*}
 for $n\geq 1$ and
$k\geq 0$, where $T_{0,0}=1$, $T_{0,k}=T_{k,-1}=0$ for $k>0$. If
$$g_kg_{k+1}- h_kf_{k+1}\geq 0$$ for all $k \geq 0$, then its
row generating functions $T_n(q)$ form a strongly $q$-log-convex
sequence.
\end{thm}

Let $[T_{n,k}]_{n,k}$ be a triangle. Assume that $m\geq n$. For
$0\le t\le m+n$, define
\begin{equation*}
   T_k(m,n,t)=T_{n-1,k}T_{m+1,t-k}+T_{m+1,k}T_{n-1,t-k}-T_{m,k}T_{n,t-k}-T_{n,k}T_{m,t-k}
\end{equation*}
if $0\le k<t/2$, and
\begin{equation*}
    T_k(m,n,t)=T_{n-1,k}T_{m+1,k}-T_{n,k}T_{m,k}
\end{equation*}
if $t$ is even and $k=t/2$. In \cite[Theorem 2.1]{ZS15}, we proved
 the next result.

\begin{thm} The transformation $y_n=\sum_{k\geq0}T_{n,k}x_k(q)$ preservers
the strong $q$-log-convexity if the following two conditions hold:
\begin{itemize}
  \item [\rm(C1)] Its row generating functions form a strongly
$q$-log-convex sequence;
  \item [\rm(C2)] There exists an index $r=r(m,n,t)$
such that $T_k(m,n,t)\ge
  0$ for $k\le r$ and $T_k(m,n,t)<0$ for $k>r$.
\end{itemize}
\end{thm}

Thus our Example \ref{example} also indicates that the C2 condition
 is necessary.

\section{Acknowledgements}
The author would like to thank the anonymous reviewers for many
valuable remarks and suggestions to improve the original manuscript.
In addition, this paper is a revision version of arXiv:1609.01544.
The results of this paper have been reported on the Fifth National
Conference on Theory of Combinatorial Numbers (Sep. 18--20, 2015,
Dalian University of Technology, Dalian) and the Workshop on
unimodality properties of Combinatorial sequences (Nov. 27--29,
2015, Nankai University, Tianjin), and Institute of Mathematics of
Academia Sinica, Taipei (Jan 19, 2016).



\end{document}